\documentclass[a4paper]{article}
\usepackage{amsmath}
\usepackage{amssymb}
\usepackage{amscd,latexsym}
\usepackage[latin1]{inputenc}
\usepackage[francais]{babel} 
\usepackage{fancyhdr}
\usepackage{amsthm}

\newtheorem{lem}{Lemme}[section]
\newtheorem{cor}[lem]{Corollaire}
\newtheorem{theo}[lem]{Th{\'e}or{\`e}me}
\newtheorem{prop}[lem]{Proposition}
\newtheorem{defi}[lem]{D{\'e}finition}
\newtheorem{rem}[lem]{Remarque}
\newtheorem*{lembis}{Lemme}
\newenvironment{dem}{\bf Démonstration : \rm}{\hfill
$\blacksquare$\par\medbreak} 
\DeclareMathOperator{\hes}{\rm{Hess}}
\DeclareMathOperator{\om}{\Omega}
\DeclareMathOperator{\omh}{\widehat{\Omega}}
\DeclareMathOperator{\h}{\widetilde{H}^1_0(\om \setminus A)}
\DeclareMathOperator{\oma}{\om \setminus A}

\DeclareMathOperator{\H1}{H^1_0(\om)}
\DeclareMathOperator{\ep}{\varepsilon}
\DeclareMathOperator{\ca}{ \rm{ Cap }}
\DeclareMathOperator{\moi}{\frac{1}{2}}

\newcommand{\pa}[1]{\frac{\partial}{\partial {#1}}}
\DeclareMathOperator{\im}{\int_\Omega}
\title{Capacité et inégalité de Faber-Krahn dans $\mathbb{R}^n$} 
\author{
J. Bertrand
\and
B. Colbois}
\date{}

\begin{document}

\maketitle            

\begin{center}{ \bf Abstract}
\end{center}
\begin{quotation}
In this paper, we define a new capacity which allows us to control the behaviour of the Dirichlet spectrum of a compact Riemannian manifold with boundary, with "small" subsets (which may intersect the boundary) removed. This result generalizes a classical result of Rauch and Taylor ("the crushed ice theorem"). In the second part, we show that the Dirichlet spectrum of a sequence of bounded Euclidean domains converges to the spectrum of a ball with the same volume, if the first eigenvalue of these domains converges to the first eigenvalue of a ball. 
\end{quotation}

\begin{center}{ \bf R\'esum\'e}
\end{center}
\begin{quotation}
Dans cet article, nous obtenons un contrôle de la variation du spectre du laplacien avec condition de Dirichlet lorsque l'on modifie une variété compacte à bord par excision d'un sous-ensemble qui intersecte éventuellement le bord. Ce résultat généralise un résultat classique de Rauch et Taylor (théorème de la glace pilée). Dans une seconde partie, nous établissons la convergence du spectre de Dirichlet d'une famille de domaines euclidiens bornés vers le spectre d'une boule de même volume, lorsque la première valeur propre de la famille converge vers celle de la boule.
\end{quotation}
\bigskip
\section{Introduction} 
Dans cet article nous nous intéressons au comportement du spectre du laplacien avec condition de Dirichlet, lorsque l'on modifie la topologie d'une variété compacte à bord par excision d'une partie. Nous nous intéressons également à la stabilité de l'inégalité de Faber-Krahn dans l'espace euclidien. Dans \cite{rauch75}, J. Rauch et M. Taylor ont démontré que le spectre de Dirichlet d'un domaine borné de l'espace euclidien, excisé par une partie compacte (contenu dans l'intérieur du domaine), converge vers le spectre de Dirichlet du domaine initial lorsque la capacité électrostatique du compact converge vers zéro (problème de la glace pilée). Il existe des versions plus explicites de ce résultat, où l'on contrôle (de manière non uniforme et par des constantes non explicites) la variation des valeurs propres en fonction de la capacité électrostatique de la partie excisée (si cette capacité est suffisamment petite), voir par exemple \cite[Theorem 6.1]{mcgil}. Nous renvoyons à \cite{noll00,daners03,flucher95} pour d'autres résultats en relation avec ce problème. Remarquons cependant que dans tous les résultats cités ci-dessus, les parties excisées sont supposées être contenues dans un compact fixé à priori, contenu dans l'intérieur du domaine. On ne peut s'affranchir de cette hypothèse puisqu'il est facile de voir que la convergence vers zéro de la capacité électrostatique n'est pas une condition nécessaire pour obtenir la convergence du spectre de la variété excisée vers celui de la variété à bord de départ. Pour s'en convaincre, il suffit de considérer un disque du plan euclidien que l'on excise par un anneau concentrique d'épaisseur petite. Lorsque l'on fait tendre l'anneau vers le bord du disque et l'épaisseur de l'anneau vers zéro, la première valeur propre de Dirichlet converge trivialement vers celle du disque alors que la capacité de l'anneau tend vers l'infini. Rappelons que la capacité électrostatique modélise la capacité d'un condensateur, elle est définie comme le minimum de l'énergie des fonctions dans $H^1_0$ qui sont égales à $1$ sur la partie excisée.

Dans la première partie de cet article, nous définissons une nouvelle capacité (au sens de la définition de G. Choquet \cite{choquet}) que nous appellerons \og capacité de Dirichlet \fg\, et nous parvenons à contrôler dans le théorème \ref{capacite} la variation du spectre de Dirichlet en fonction de la capacité de Dirichlet et de constantes explicites. De plus, la proximité des spectres implique que la capacité de Dirichlet du sous-ensemble excisé est petite. Cette nouvelle capacité nous permet de traiter le cas où le sous-ensemble excisé intersecte le bord de la variété considérée.

Dans la deuxième partie de cet article, nous appliquons le résultat obtenu dans la première partie pour démontrer un résultat de stabilité de l'inégalité de Faber-Krahn dans un espace euclidien de dimension au moins 3. L'inégalité de Faber-Krahn établit que, pour des domaines bornés à bord lisse de volume fixé, la première valeur propre de Dirichlet est minimale lorsque le domaine est une boule euclidienne et que ce minimum est caractéristique de la boule. Sous les hypothèses de l'inégalité de Faber-Krahn, nous montrons dans le théorème \ref{convergence} que l'on peut contrôler de manière explicite toute partie finie du spectre de Dirichlet d'un domaine en fonction des valeurs propres correspondantes de la boule de même volume, si la première valeur propre du domaine est suffisamment proche de celle de la boule. En conséquence, il y a convergence du spectre de Dirichlet d'une famille de domaines vers celui de la boule lorsque la première valeur propre de cette famille converge vers celle de la boule. Ce résultat de convergence a été obtenu précédemment en dimension 2 par F. Crevoisier, pour des domaines dont la topologie est contrôlée (précisément il suppose que les groupes fondamentaux des domaines considérés admettent des systèmes de générateurs dont le cardinal est borné par une constante $N$ donnée à priori) \cite{crevoisier}. L'utilisation de la capacité de Dirichlet permet de traiter le cas général. La démonstration du théorème \ref{convergence} repose sur une inégalité isopérimétrique quantitative dans l'espace euclidien établie par R.R. Hall \cite{hall} (voir le théorème \ref{hal} pour un énoncé), sur un résultat de transplantation de fonctions propres (lemme \ref{a4l3}) et sur le théorème \ref{capacite} de la première partie. Notons que sous les hypothèses du théorème \ref{convergence}, il n'y a pas de résultat de stabilité pour la distance de Hausdorff, on ne peut en déduire qu'une description géométrique grossière des domaines vérifiant ces hypothèses (voir la proposition \ref{pov} pour un énoncé, \cite{povel,sznitman97}) mais cela n'implique, à priori, aucun résultat sur le spectre de Dirichlet du domaine. Signalons également que l'inégalité de Faber-Krahn est vérifiée sur les autres espaces modèles de courbure constante. La preuve du théorème \ref{convergence} peut-être adaptée à ces autres cas, si l'on démontre une inégalité isopérimétrique analogue à celle de R.R. Hall sur ces espaces.

\section{Une capacité pour le problème de Dirichlet}
L'objet de cette partie est de définir une capacité qui permette de contrôler la variation du spectre de Dirichlet d'une variété compacte connexe à bord lisse lorsque que l'on modifie la variété par excision d'une partie. La définition de cette capacité est inspiré de la capacité électrostatique, habituellement utilisée dans ce type de problème et d'une capacité introduite par G. Courtois \cite{courtois1995} pour traiter le cas des variétés compactes sans bord.

\begin{defi}[Capacité de Dirichlet]\label{cap}Soit $\Omega$ une variété connexe, compacte, à bord $\partial \Omega$ lisse et $A$ une partie de $\Omega$. On note $\h$ l'adhérence dans $H^1_0(\om)$ de $\{ f \in C^{\infty}_c(\om \setminus \partial \Omega); f=0 \mbox{ sur un voisinage ouvert de } A\}$ et $\phi_1$  la première fonction propre du laplacien sur $\om$ avec condition de Dirichlet sur le bord, positive et unitaire pour la norme $L^2$ sur $\om$. On définit la capacité de Dirichlet de $A$ par
$$ \ca (A) = \inf \left\{ \int_{\om} |\nabla f|^2, f-\phi_1 \in \h \right\}.$$
\end{defi}

Dans cette partie nous démontrons le théorème suivant (voir \cite{courtois1995} pour un théorème analogue dans le cas sans bord).
\begin{theo}\label{capacite} Soit $\om$ une variété compacte, connexe, à bord lisse. Il existe des   
 constantes positives $B_1,(C_k)_{k \geq 1},(\ep_k)_{k \geq 1}$ telles que pour toute partie $A$ de $\Omega$, on a les inégalités
 $$
   \ca (A) \leq B_1\big(\lambda_1(\oma)- \lambda_1(\om)\big)  $$
 et si $\ca (A) < \ep_k$, 
 $$  0 \leq \lambda_k (\oma) - \lambda_k (\om) \leq C_k \ca^{\frac{1}{2}} (A).  $$
Les constantes sont définies comme suit\\
\begin{center}
$ \displaystyle{B_1 = \frac{\lambda_1(\om)+ \lambda_2(\om)}{\lambda_2(\om)-\lambda_1(\om)}}$,\\
\smallskip
$ \displaystyle{\ep_1 =\frac{\lambda_1(\om)}{16}}$, \\
\smallskip
$C_1 = 14 \sqrt{\lambda_1(\om)}$. \\
\end{center}
Pour tout entier $k \geq 2$,

$$\ep_k = \frac{\lambda_1^2(\om)}{4k^2(m_k+\sqrt{\lambda_1(\om)})^4},$$
\begin{equation}\label{a4e18}
 C_k=
2k\left(2\lambda_k(\om)\Big(1+ \frac{m_k}{\sqrt{\lambda_1(\om)}}\Big)^2+
\Big(m_k+ \frac{M_k}{\sqrt{\lambda_1(\om)}}+\sqrt{\lambda_k(\om)}\Big)^2\right)
\end{equation}
où $M_k=M_k(\om)= \max_{1 \leq i \leq k} \big|\big|\nabla \big(\frac{\phi_i}{\phi_1}\big)\big|\big|_{L^{\infty}}$, $m_k=m_k(\om)= \max_{1 \leq i \leq k} ||\frac{\phi_i}{\phi_1}||_{L^{\infty}}$ et $(\phi_i)_{i \geq 1}$ est une base orthogonale de fonctions propres de $\H1$, unitaires pour la norme $L^2$ sur $\om$.
\end{theo}

Le théorème \ref{capacite} est démontré dans le paragraphe \ref{dem_capa}. Nous montrerons en particulier que les quantités $m_k(\om)$ et $M_k(\om)$ sont bien définies. Dans le prochain paragraphe, nous établissons les propriétés de la capacité de Dirichlet nécessaires pour démontrer ce résultat.

\subsection{Premières propriétés de la capacité de Dirichlet}    
\begin{prop}\label{a3p1}Soit $\om$ une variété connexe, compacte, à bord lisse et $A$ une partie de $\om$. Il existe une unique fonction $f_A$ telle que
$$ \ca (A) = \int_{\om} |\nabla f_A|^2, \mbox{ avec } f_A-\phi_1 \in \h. $$
De plus, les conditions suivantes sont équivalentes 
$$
\begin{array}{rl}

i) & \ca (A)= 0 \\
ii) & f_A = 0 \\
iii)& \h = H^1_0 (\om)
\end{array}$$
\end{prop}
\begin{dem}
Notons $(\cdot,\cdot)$ le produit scalaire sur $\H1$ défini par 
$$(f,h) = \int_{\om} \langle\nabla f, \nabla h\rangle , $$
pour $f,h$ appartenant à $\H1$. Muni de ce produit scalaire, $\H1$ est un espace de Hilbert et pour toute partie $A$ de $\om$, $\ca (A)$ est le carré de la distance de la fonction identiquement nulle au sous-espace affine (fermé) $\phi_1 + \h$. Par conséquent la quantité $\ca (A)$ est réalisée par une unique fonction appartenant à $\H1$, que nous noterons $f_A$. On en déduit également l'équivalence des points i) et ii). Supposons maintenant que $f_A=0$, ceci implique que $\phi_1$ appartient à $\h$. Par conséquent, il existe une suite $(v_n)_{n \in \mathbb{N}}$ de fonctions lisses de $\H1$ nulles sur un voisinage ouvert de $A$, qui tend vers la fonction $\phi_1$ pour la norme usuelle de $\H1$. D'après le théorème nodal de Courant, $\phi_1$ ne s'annule pas sur $\om$ (voir \cite{chavel} pour une démonstration).  Soit $f$ une fonction lisse à support compact dans $\om\setminus \partial \om$, la suite $(\frac{f}{\phi_1}v_n)_{n \in \mathbb{N}}$ tend vers $f$ dans $\H1$ et appartient à $\h$, ce qui termine la démonstration de l'implication ii) donne iii) par densité. L'implication réciproque est immédiate puisque l'hypothèse iii) implique que $\phi_1$ appartient à $\h$.

\end{dem}

La capacité de Dirichlet est une capacité au sens de la définition de G. Choquet \cite{choquet}. Nous démontrons cette propriété dans la proposition ci-dessous.

\begin{prop} Soit $\om$ une variété connexe, compacte, à bord lisse. La capacité de Dirichlet vérifie les propriétés suivantes :

\begin{tabular}{cl}
1) & $ \ca(A) \leq \int_{\om} |\nabla \phi_1|^2< +\infty$ pour toute partie $A$ de $\om$.\\
2) & $\ca(A) \leq \ca(B)$ pour toutes parties  $A,B$ de $\om$ vérifiant $A \subset B$. \\
3) & $\ca(A)= \lim_{n \rightarrow +\infty} \ca(A_n)$ pour tout partie $A$ vérifiant \\ & $A=\cup_n A_n$ 
  avec  $(A_n)_{n \in \mathbb{N}}$ une suite croissante de parties de $\om$. \\
4) & $\ca(K) = \lim_{n \rightarrow +\infty} \ca(K_n)$ pour toute partie compacte $K$ vérifiant \\
 & $K=\cap_n K_n$ avec  $(K_n)_{n \in \mathbb{N}}$ une suite décroissante de parties compac\-tes \\ 
 & de $\om$. \\
\end{tabular}
\end{prop}
\begin{dem}
La démonstration qui suit est inspirée de la preuve de la proposition 2.4 de \cite{courtois1995} qui traite le cas des variétés compactes sans bord.

Les propriétés 1) et 2) découlent de la définition de la capacité de Dirichlet. Notons $A=\cup_n A_n$ vérifiant les hypothèses du point 3) et $F_{A}= \phi_1 + \h$. Par définition, $\ca(A_n)$ est une suite croissante majorée par $\ca(A)< +\infty$, par conséquent elle converge. Nous allons en déduire la convergence de la suite de fonctions $(f_{A_n})_{n\geq 1}$. Pour cela, remarquons que $(F_{A_n})_{n \geq 1}$ est une suite décroissan\-te pour l'inclusion, par conséquent pour tout entier $n$ et tout entier $m \geq n$, $f_{A_n}-f_{A_m}$ appartient à $\widetilde{H}^1_0(\om \setminus A_n)$. Or $f_{A_n}$ est le projeté orthogonal pour le produit scalaire $(\cdot,\cdot)$, de $0$ sur $F_{A_n}$, par conséquent le théorème de Pythagore et la convergence de la suite $(\ca (A_n))_{n \geq 1}$ impliquent la convergence de la suite $(f_{A_n})_{n \geq 1}$ vers une fonction $u$. Par construction, la fonction $u-\phi_1$ appartient à $ \cap_n \widetilde{H}^1_0(\om \setminus A_n)$. La propriété 3) découle alors de l'égalité 
$$   \cap_n \widetilde{H}^1_0(\om \setminus A_n)= \h.$$ 
La démonstration de cette égalité est similaire à celle donnée par G. Courtois dans le cas sans bord, nous renvoyons donc à \cite[proposition 2.4]{courtois1995} pour une preuve.

Soit $K=\cap_n K_n$ un compact de $\om$ vérifiant les hypothèses de 4). D'après la propriété 2), $\ca (K_n) \geq \ca(K)$ pour tout entier $n$. Fixons un réel $\ep$ positif et
$f$ une fonction lisse de $\H1$, nulle sur un voisinage ouvert $V$ de $K$, vérifiant 
$$ \int_{\om} |\nabla (\phi_1+ f)|^2 - \ca(K) \leq \ep.$$
Pour tout entier $n$ assez grand, $ K_n$ est contenu dans $V$ et par définition de la capacité,
$$ \ca (K_n) \leq      \int_{\om} |\nabla (\phi_1+ f)|^2,$$
ce qui termine la démonstration.  

\end{dem}

La capacité de Dirichlet vérifie également une propriété de sous-additivité.
\begin{lem}\label{a4l6}
Soit $A$ un sous-ensemble de $\om$, $f_A$ la fonction définie dans la proposition \ref{a3p1} et qui vérifie $\ca(A)=\int_{\om} |\nabla f_A|^2$. Sous ces hypothèses, on a l'inégalité 
$$ f_A \leq \phi_1,$$
où $\phi_1$ est la fonction propre associée à la première valeur propre de $\om$ introduite dans la définition \ref{cap}.

Soit $A,B$ deux parties de $\om$. On a l'inégalité
$$\ca (A\cup B) \leq \ca(A) + \ca(B).$$
\end{lem} 

\begin{dem}
Démontrons la première partie de l'énoncé. Soit $f_A$ la fonction dont l'énergie réalise $\ca (A)$. Notons $h= \min (f_A,\phi_1)$. Par construction, la fonction $h$ appartient à $\H1$ (voir \cite[lemme 7.6]{GT}) et vérifie 
$$ \int_{\om} |\nabla h|^2 \geq \ca (A).$$
Par conséquent, on en déduit
\begin{equation}\label{truc}
 \int_{\om \cap \{\phi_1-f_A \geq 0\}^c}|\nabla f_A|^2 \leq \int_{\om \cap \{\phi_1-f_A \geq 0\}^c}|\nabla \phi_1|^2.
\end{equation}
Considérons maintenant la fonction $\psi= \max (f_A,\phi_1)$. Par construction, $\int_{\om} \psi^2 \geq \int_{\om}\phi_1^2$ et en utilisant (\ref{truc}), on en déduit
$$ \int_{\om}|\nabla \psi|^2 \leq \int_{\om}|\nabla \phi_1|^2.$$
Par conséquent, la caractérisation variationnelle de la première valeur propre entraîne l'existence d'un réel $\mu \geq 1$ tel que $\psi=\mu \phi_1$. En utilisant la définition de $f_A$, on montre que $\mu=1$; ce qui démontre la première inégalité.

Soit $f_B$ la fonction dont l'énergie réalise $\ca (B)$. La première partie du lemme montre que la fonction $\theta = \max (f_A,f_B)$ est une fonction test pour la capacité de la partie $A \cup B$. De plus, pour presque tout $x$ dans $\om$, le gradient de $\theta$ vérifie 
$$ |\nabla \theta|^2 \leq |\nabla f_A|^2 + |\nabla f_B|^2,$$
ce qui permet de conclure.
\end{dem}

\subsection{Démonstration du théorème \ref{capacite}}\label{dem_capa}
Démontrons la première partie du théorème. On note $\psi = \phi_1 -\phi_A$ où $ \phi_A$ est la première fonction propre 
de Dirichlet du domaine $\oma$, positive et unitaire pour la norme $L^2$. Précisons que si le sous-ensemble $A$ n'est pas compact, le spectre de Dirichlet de $\oma$ désigne le spectre de la forme quadratique usuelle sur l'espace $\h$. Par choix de $\psi$, on a l'inégalité 
$$ \ca(A) \leq \im |\nabla \psi|^2.$$
En développant puis en utilisant la formule de Green, on obtient
\begin{equation}\label{a4e19}
\ca(A) \leq \lambda_1(\om) + \lambda_1(\oma) -2 \lambda_1(\om)\im \phi_1\phi_A.
\end{equation}

Notons $(c_i)_{i \geq 1}$ les coefficients du développement en série de Fourier de la fonction $\phi_A$ dans la base de fonctions propres $(\phi_k)_{k\geq 1}$ de $\H1$. Avec ces notations, on a $c_1 = \im \phi_1\phi_A>0$. Par choix de la normalisation de $\phi_A$, on a $||(c_i)||_{l^2(\mathbb{N}^*)}=1$, on en déduit
$$ \lambda_1(\oma) \geq \lambda_1(\om)c_1^2 + \lambda_2(\om)(1-c_1^2).$$
En utilisant $0<c_1\leq 1$ , il vient
$$ 1-c_1 \leq  \frac{\lambda_1(\oma) -\lambda_1(\om)}{\lambda_2(\om)-\lambda_1(\om)}.$$
En injectant cette inégalité dans (\ref{a4e19}), on obtient 
$$ \ca(A) \leq \left(\lambda_1(\oma)- \lambda_1(\om)\right)\left(1 + \frac{2\lambda_1(\om)}{\lambda_2(\om)-\lambda_1(\om)}\right).$$

Montrons maintenant la deuxième partie de l'énoncé. Fixons un entier $k \geq 2$. L'inégalité de gauche est immédiate. Nous montrons l'inégalité de droite en appliquant le principe du min-max au spectre de $\oma$. On note $E_k$ le sous-espace vectoriel de $\h$, définie par
$$ E_k = \left\{ f(1- \frac{f_A}{\phi_1}), f \in \mbox{ \rm{Vect}}(\phi_1,\cdots,\phi_k) \right\}.$$
Le fait que l'ensemble $E_k$ soit contenu dans $\h$ est une conséquence du lemme ci-dessous.  
\begin{lem}\label{lem} Sous les hypothèses du théorème \ref{capacite} et en notant $\phi_k$ la $k^e$ fonction propre de Dirichlet sur $\Omega$, nous avons les estimations 
$$ \sup_{\Omega} \frac{|\phi_k|}{\phi_1} < \infty \mbox{ et } \sup_{\Omega} \left|\nabla \left(\frac{\phi_k}{\phi_1}\right) \right| < \infty. $$
\end{lem}
\noindent Nous donnons une démonstration de ce lemme à la fin de cette section.

\smallskip

Notons $\psi_k= \phi_k(1-\frac{f_A}{\phi_1})$ et pour alléger les notations, \\$M_k=M_k(\om)= \max_{1 \leq i \leq k} \big|\big|\nabla \big(\frac{\phi_i}{\phi_1}\big)\big|\big|_{L^{\infty}}$ et $m_k=m_k(\om)= \max_{1 \leq i \leq k} ||\frac{\phi_i}{\phi_1}||_{L^{\infty}}$. Pour simplifier, nous noterons $(\lambda_i)_{i \geq 1}$ les valeurs propres de $\om$. Soit $i,j$ deux entiers dans $\{1,\cdots,k\}$, on a l'égalité 
$$ \im \psi_i\psi_j= \delta_{i,j} + \im \frac{\phi_i }{\phi_1}\frac{\phi_j }{\phi_1}f_A^2- 2 \im \frac{\phi_i\phi_j}{\phi_1}f_A.$$
En utilisant l'inégalité de Poincaré, on obtient
$$
\left| \im \psi_i\psi_j- \delta_{i,j}\right| \leq \frac{m_k^2}{\lambda_1}\ca (A) +2 \frac{m_k}{\sqrt{\lambda_1}}\ca^{\frac{1}{2}}(A).
$$
Quitte à supposer $\ca (A)\leq 1$, il vient
\begin{equation}\label{a4e1}
\left| \im \psi_i\psi_j- \delta_{i,j}\right| \leq \left(1+ \frac{m_k}{\sqrt{\lambda_1}}\right)^2\ca^{\frac{1}{2}}(A).
\end{equation}
Notons   $X_1= \left(1+ \frac{m_k}{\sqrt{\lambda_1}}\right)^2\ca^{\frac{1}{2}}(A).$

\begin{multline}
\im \langle\nabla \psi_i,\nabla \psi_j\rangle = \lambda_i \delta_{ij}- \overbrace{\im f_A \left(\langle\nabla (\frac{\phi_i}{\phi_1}), \nabla \phi_j\rangle  + \langle\nabla (\frac{\phi_j}{\phi_1}), \nabla \phi_i\rangle  \right)}^{(I)} \\
+ \overbrace{\im f_A^2 \langle\nabla (\frac{\phi_i}{\phi_1}),\nabla (\frac{\phi_j}{\phi_1})\rangle }^{(II)}
- \overbrace{\im \langle\nabla f_A, \frac{\phi_i}{\phi_1}\nabla \phi_j+ \frac{\phi_j}{\phi_1}\nabla \phi_i\rangle }^{(III)}\\
+ \overbrace{\im \frac{\phi_i\phi_j}{\phi_1^2}|\nabla f_A|^2}^{(IV)} + \overbrace{\im \langle f_A\nabla f_A, \frac{\phi_i}{\phi_1}\nabla \left( \frac{\phi_j}{\phi_1}\right)+ \frac{\phi_j}{\phi_1}\nabla \left( \frac{\phi_i}{\phi_1}\right)\rangle }^{(V)}.
\end{multline}

On majore chacun des termes comme suit.
$$ |(I)| \leq 2M_k\left(\frac{\lambda_k}{\lambda_1}\right)^{\moi}\ca^{\moi}(A) \mbox{, } |(II)| \leq \frac{M_k^2}{\lambda_1}\ca(A),$$
$$ |(III)| \leq 2m_k \lambda_k^{\frac{1}{2}}\ca^{\frac{1}{2}}(A) \mbox{, } |(IV)| \leq m_k^2 \ca(A),$$
$$ |(V)| \leq \frac{2m_kM_k}{\lambda_1^{\frac{1}{2}}}\ca(A).$$   
Ce qui donne 
$$
\left|\im \langle\nabla \psi_i,\nabla \psi_j\rangle - \lambda_i \delta_{ij}\right| \leq 2\left(\left(\frac{\lambda_k}{\lambda_1}\right)^{\moi}M_k+m_k\sqrt{\lambda_k}\right)\ca^{\moi}(A) + \left(\frac{M_k}{\sqrt{\lambda_1}}+ m_k\right)^2 \ca(A).
$$
On majore cette expression par
\begin{equation}\label{a4e2}
\left|\im \langle\nabla \psi_i,\nabla \psi_j\rangle - \lambda_i \delta_{ij}\right| \leq \left(m_k+ \frac{M_k}{\sqrt{\lambda_1}}+\sqrt{\lambda_k}\right)^2\ca^{\moi}(A).
\end{equation} 
Notons   $ X_2=\left(m_k+ \frac{M_k}{\sqrt{\lambda_1}}+\sqrt{\lambda_k}\right)^2\ca^{\moi}(A).$

Il reste à établir une condition suffisante pour que la famille des fonctions tests $(\psi_i)_{1 \leq i \leq k}$ soit de dimension $k$. Supposons 
$$\left| \im \psi_i\psi_j- \delta_{i,j}\right| \leq X$$
et $\sum_{i=1}^k a_i \psi_i=0,$ avec $(a_i)_{1 \leq i \leq k}$ une famille de réels vérifiant $\sum_{i=1}^k a_i^2 =1$.\\
Nous pouvons supposer $ |a_1|^2 \geq \frac{1}{k}$. On déduit de l'égalité $ \im (\sum_{i=1}^ka_i\psi_i)\psi_1=0$, l'estimation 
$$\int_{\om}\psi_1^2 \;|a_1| \leq \sum_{i=2}^k |a_i|\left|\im \psi_i\psi_1  \right|.$$
On en déduit $ X \geq \frac{1}{k}.$ Par conséquent, en utilisant (\ref{a4e1}), nous en déduisons que la famille $(\psi_i)_{1\leq i \leq k}$ est libre si la capacité de l'ensemble $A$ vérifie
\begin{equation}\label{a4e17}
  \ca(A) < \frac{\lambda_1^2}{4k^2(m_k+\sqrt{\lambda_1})^4},
\end{equation}  
ce que l'on suppose dorénavant.\\
Soit $f$ appartenant à $E_k$. Supposons que $f= \sum_{i=1}^k b_i\psi_i$ avec $\sum_{i=1}^k b_i^2=1$. On déduit des estimations précédentes, les inégalités
$$ \int_{\oma} |\nabla f|^2 \leq \lambda_k(\om) + kX_2,$$
$$ \int_{\oma} f^2 \geq 1 -kX_1.$$
Par conséquent en utilisant (\ref{a4e17}), on en déduit
$$ \lambda_k(\oma) \leq \lambda_k(\om) + 4k\lambda_k(\om)X_1+ 2kX_2,$$ 
ce qui donne le résultat. 
La démonstration dans le cas $k=1$ est similaire, elle utilise la fonction test $\phi_1-f_A$, nous laissons les détails au lecteur.

\smallskip

\noindent
{\bf Démonstration du lemme \ref{lem}}

Notons $h_k = \frac{\phi_k}{\phi_1}$. D'après le théorème nodal de Courant, le supremum de $h_k$ est fini sur toute partie compacte de $\om\setminus \partial \Omega$. D'autre part par le principe du maximum \cite{chavel}, le gradient de $\phi_1$ vérifie en tout point $x$ du bord $\partial \om$, 
\begin{equation}\label{a3e2}
 \eta(x).  \phi_1 >0
\end{equation}
où $\eta(x)$ est la normale unitaire intérieure en $x$. Pour contrôler $h_k$, nous utilisons des coordonnées normales par rapport à $\partial \om$. Par abus de notation, on note $\phi_k(t,x)$ la fonction $\phi_k$ dans ces coordonnées (avec $t$ la coordonnée normale et $x$ appartient à $\partial \om$). Sous les hypothèses du lemme, toute fonction propre du laplacien pour les conditions de Dirichlet au bord, est lisse sur $\om$ et nulle sur $\partial \om$. Par conséquent pour tout entier positif $k$ et tout $x$ dans $\partial \om$ fixés, un développement limité donne
\begin{equation}\label{a3e1}
 \phi_k(t,x) = t \pa{t}\phi_k(0,x) + O (t^2),    
\end{equation}
où $O(t^2)$ dépend de $||\hes (\phi_k)||_{L^{\infty}}$, on en déduit donc
$$ \sup_{\Omega} \frac{|\phi_k|}{\phi_1} < \infty.$$

Estimons maintenant le gradient de $h_k$. Un calcul donne
$$ \nabla h_k = \frac{1}{\phi_1^2}(\phi_1 \nabla \phi_k - \phi_k \nabla \phi_1).$$
Donc, il suffit également de contrôler $|\nabla h_k|$ au voisinage de $\partial \om$. Dans ce but, nous utilisons la même méthode que ci-dessus.   
\begin{multline}
(\phi_1 \nabla \phi_k - \phi_k \nabla \phi_1)(t,x) = (\phi_1\partial_t\phi_k -\phi_k\partial_t\phi_1)(t,x)\pa{t} + \\ \sum_{i=1}^{n-1} (\phi_1\partial_{x_i}\phi_k - \phi_k\partial_{x_i}\phi_1)(t,x)\pa{x_i}.
\end{multline} 
D'après (\ref{a3e1}), on en déduit 
$$ (\phi_1\partial_t\phi_k -\phi_k\partial_t\phi_1)(t,x)= O(t^2),$$
où $O(t^2)$ dépend de $ ||\hes (\phi_k)||_{L^{\infty}}$ et de $||\hes (\phi_1)||_{L^{\infty}}$. On obtient une formule similaire pour les autres termes, par conséquent $\sup_{\Omega} \left|\nabla \left(\frac{\phi_k}{\phi_1}\right) \right| < \infty$.

\section{Inégalité de Faber-Krahn}
L'inégalité de Faber-Krahn affirme que parmi les domaines bornés de l'espace euclidien de volume fixé et dont le bord est lisse, la première valeur propre de Dirichlet est minimale lorsque le domaine est une boule et que ce minimum est caractéristique de la boule. A. Melas a démontré qu'un domaine convexe dont la première valeur propre de Dirichlet est proche de celle d'une boule de même volume, est proche au sens de la distance de Hausdorff, d'une boule de même volume. 
\begin{theo}[Melas, \cite{melas}]\label{mel}Soit $\Omega$ un domaine convexe  borné de l'espace euclidien dont le bord est lisse.
Supposons que la première valeur propre de Dirichlet de ce domaine vérifie
$$ \lambda_1 (\Omega) \leq (1+\ep)\lambda_1(B)$$
où $B$ est une boule euclidienne de même volume que $\Omega$. Sous ces hypothèses, il existe une fonction $c(\ep)$ explicite, dépendant de la dimension et tendant vers $0$ avec $\ep$, telle que la distance de Hausdorff entre $\Omega$ et $B$ vérifie
$$ d_H(\Omega,B) \leq c(\ep)vol (\om)^{\frac{1}{n}}.$$
\end{theo}

En utilisant la décroissance du spectre de Dirichlet par rapport à l'inclusion, on déduit de ce résultat que sous les hypothèses du théorème de A. Melas, tout le spectre du convexe est (non uniformément) proche de celui d'une boule de même volume.
\begin{cor}Sous les hypothèses du théorème \ref{mel}, il existe une fonction explicite $r(\ep)$ telle que pour tout entier $k\geq 2$, on a l'estimation
$$(1-r(\ep))\lambda_k(B)\leq \lambda_k(\om)\leq (1+r(\ep))\lambda_k(B).$$
\end{cor}
L'objet de cette partie est d'établir un résultat similaire au corollaire ci-dessus sans hypothèse de convexité. Signalons que sans cette hypothèse de convexité, on ne peut espérer obtenir ce type de résultat comme un corollaire d'un résultat de proximité géométrique (au sens de Hausdorff). Précisément, nous allons montrer le
\begin{theo}\label{convergence}Soit $\Omega$ un domaine  borné de l'espace euclidien de dimension $n$ ($n \geq 3$) dont le bord est lisse.
Supposons que la première valeur propre de Dirichlet de ce domaine vérifie
$$ \lambda_1 (\Omega) \leq (1+\ep)\lambda_1(B)$$
où $B$ est une boule euclidienne de même volume que $\Omega$. Sous ces hypothèses, il existe des constantes positives $(\alpha_k)_{k\geq2}$ et $(C'_k)_{k \geq 2}$ telles que si $\ep < \alpha_k$ alors
$$|\lambda_k(B)-\lambda_k(\om)|\leq C'_k \ep^{\frac{1}{80n}},$$
où les constantes $\alpha_k$ et $C'_k$ dépendent des paramètres $n,k,vol (\om)$.
\end{theo}  
\begin{rem}La démonstration du théorème \ref{convergence} permet d'obtenir des expressions explicites des constantes $C'_k$ et $\alpha_k$ ci-dessus, dépendant des paramètres $n,k,vol (\om)$ et des quantités $(\lambda_i(B))_{i \geq 1}, m_k(B),M_k(B)$. Cependant, les expressions obtenues sont trop compliquées pour être reproduites ici.
\end{rem}
\begin{rem}La démonstration du théorème \ref{convergence} permet également de traiter le cas des domaines plans. Il suffit d'utiliser une inégalité isopérimétrique de Bonnesen à la place de celle de Hall. La démonstration ci-dessous peut aussi être adaptée aux cas de la sphère et de l'espace hyperbolique si l'on démontre sur ces espaces une inégalité isopérimétrique analogue au théorème \ref{hal}.
\end{rem}
La démonstration du théorème \ref{convergence} fait l'objet des paragraphes suivants.

\subsection{Inégalité isopérimétrique de Hall}\label{hall}
La première étape de la démonstration du théorème \ref{convergence} consiste à utiliser une inégalité isopérimétrique quantitative établie par R.R. Hall \cite{hall}.
\begin{theo}[{\cite[theorem 1]{hall}}]\label{hal} Soit $\om$ un domaine borné de l'espace euclidien de dimension $n$ ($n \geq 3$), dont le bord est lisse. On note $F(\om)$ la quantité définie par
$$ \sup_{\cal A}\big(vol (B\cap \om)\big)=\big(1-F(\om)\big)vol(\om),$$
où ${\cal A}=\{B \mbox{ boules de volume } vol (B)= vol (\om)\}$.

Il existe une constante $c(n)>0$ ne dépendant que de la dimension $n$ telle que
$$ vol(\partial \om) \geq vol (\partial B)(1+c(n)F(\om)^4).$$
\end{theo}
\begin{rem}Le théorème 1 de \cite{hall} comporte l'hypothèse supplémentaire $F(\om) \leq 1-\frac{1}{(5n)^2}$, R.R. Hall précise que l'on peut supprimer cette hypothèse, quitte à modifier la constante dans la conclusion du théorème.
\end{rem}

On déduit de cette inégalité isopérimétrique une description géométrique approximative des domaines euclidiens vérifiant les hypothèses du théorème \ref{convergence}.

\begin{prop}[{\cite[Theorem A]{povel}}]\label{pov}Soit $\om$ un domaine borné de l'espace euclidien de dimension $n\geq 3$, dont le bord est lisse. Il existe une constante ne dépendant que de la dimension $c(n)$ telle que, pour tout $h$ dans $]0,vol (\om)^{-\moi}[$, il existe une boule $B$ de même volume que $\om$ pour laquelle on l'inégalité
$$ vol (\om \setminus B) \leq  vol(\om)^{\frac{7}{8}}\left(c(n)\Big(\frac{\ep^{\moi}}{h}\Big)^{\frac{1}{4}}+ vol(\om)^{\frac{1}{8}}\left(1-\Big(\frac{1-h (vol(\om))^{\moi}}{(1+\ep)^{\moi}}\Big)^n\right)\right),$$
où $\ep= \frac{\lambda_1(\om)}{\lambda_1(B)}-1$ et $c(n)$ est la constante introduite dans le théorème \ref{hal}.
\end{prop} 
\begin{rem}Cette proposition montre qu'un domaine borné dont la première valeur propre est proche de celle d'une boule de même volume, est, à des ensembles de petits volumes près, une boule. Cependant, cela n'implique aucun contrôle du diamètre du domaine ni, à priori, des valeurs propres suivantes. On peut penser, par exemple, à un disque du plan auquel on colle de manière lisse un \og rectangle \fg\, long et fin. 
\end{rem}

On déduit de cette proposition une majoration de la forme 
\begin{equation}\label{a4e10}
 vol (\om \setminus B) \leq c_1(n)vol (\om) \ep^{\frac{1}{10}}.
\end{equation}
{\bf Démonstration de l'inégalité (\ref{a4e10}).}

Posons $h=\frac{\ep^{\alpha}}{vol (\om)^{\moi}}$ avec $0<\alpha<\moi$. On suppose $0<\ep<1$. En utilisant la minoration
$$\frac{1}{(1+\ep)^{\moi}} 
\geq 1-\ep,$$
on obtient par hypothèse sur $\alpha$, 
$$\frac{1-h (vol(\om))^{\moi}}{(1+\ep)^{\moi}}=\frac{1-\ep^{\alpha}}{(1+\ep)^{\moi}}\geq (1-\ep^{\alpha})^2.$$
On en déduit
$$ 1-  \left(\frac{1-h (vol(\om))^{\moi}}{(1+\ep)^{\moi}}\right)^{n} \leq 2n \ep^{\alpha},$$
d'où la majoration
$$ vol (\om \setminus B) \leq vol (\om)(c(n)+2n)\ep^{\min \{1/4(1/2-\alpha),\alpha\}}.$$
La majoration est optimale en choisissant $\alpha=\frac{1}{10}$. Notons $c_1(n)=c(n)+2n$.
\begin{flushright}
$\blacksquare$
\end{flushright}

Soit $\om$ un domaine vérifiant les hypothèses du théorème \ref{convergence}. Fixons $B_0$ une boule de même volume que $\om$ telle que $1-F(\om)= \frac{vol(B_0\cap \om)}{vol (\om)}$. Notons $D_1=B_0 \cap \om$, le domaine $B_0$ après excision par le domaine $A_0=B_0 \cap \om^c$ (où $\om^c$ est le complémentaire de $\om$ dans $\mathbb{R}^n$). Pour décrire la suite de la démonstration, nous avons besoin des estimations ci-dessous. 
 $$
 \begin{array}{ccl}
 |\lambda_k(\om)- \lambda_k(B_0)| & \leq & |\lambda_k(\om)-\lambda_k( D_1)| + |\lambda_k(D_1)-\lambda_k(B_0)| \\
  & \leq & \lambda_k( D_1) - \lambda_k(\om) + \lambda_k( D_1)-\lambda_k(B_0) \\
  & \leq & \lambda_k(B_0) - \lambda_k(\om) +2(\lambda_k(D_1)-\lambda_k(B_0)), 
\end{array}
$$
on obtient finalement
\begin{equation}\label{a4e15}
|\lambda_k(\om)- \lambda_k(B_0)| \leq  \lambda_k(B_0)- \lambda_k(\om \cup B_0)  +2(\lambda_k( D_1)-\lambda_k(B_0)).
\end{equation}
La suite de la démonstration repose sur l'idée suivante : considérons un domaine contenant une boule de volume presque égal à celui du domaine. Alors on peut montrer que le spectre du domaine est proche de celui de la boule. La méthode consiste à construire à l'aide des fonctions propres du domaine que l'on tronque par des fonctions plateaux, des fonctions tests sur une boule concentrique et de rayon légèrement supérieur à celui de la boule contenue dans le domaine. On conclut en utilisant l'expression explicite du spectre d'une boule euclidienne en fonction de son rayon. 
Cette méthode permet également de contrôler la première valeur propre d'une boule de rayon légèrement supérieur à celui de $B_0$, excisée par $A_0$, en fonction de la première valeur propre du domaine $\om$ et d'un terme petit controlé. Nous déduirons du contrôle de cette première valeur propre et du théorème \ref{capacite}, une majoration du terme $ \lambda_k( D_1)-\lambda_k(B_0)$.

\subsection{Transplantation de fonctions propres}
Dans la suite, toutes les boules considérées auront même centre que la boule $B_0$ et on choisira ce centre comme origine de $\mathbb{R}^n$. De plus, on notera $R$ le réel positif tel que 
$$ vol (B_0)= vol (B(R)) \mbox{  soit  } R= \left(\frac{vol (\om)}{w_n}\right)^{\frac{1}{n}},$$
où $w_n$ désigne le volume de la boule unité de $\mathbb{R}^n$.
Dans ce paragraphe, on établit le résultat technique suivant. 
\begin{lem}\label{a4l3}Soit $\omh$ un domaine borné de $\mathbb{R}^n$ tel que $vol (\omh \setminus B_0) \leq \theta$ et $k$ un entier non nul fixé. Il existe des constantes $R'=R'(R,\theta,n)$, $\gamma=\gamma(k,n,\lambda_k(\omh))$ et $\kappa=\kappa(k,n,\lambda_k(\omh))$ telles que la propriété suivante est vérifiée. Pour $\theta < \kappa$, il existe une famille de fonctions $(F_i)_{1 \leq i \leq k}$ de $H^1_0(B')$, orthonormée pour le produit scalaire $L^2$ sur $B'$, qui vérifie pour tout $i$ dans $\{1,\cdots,k\}$, 
$$ \int_{B'} |\nabla F_i|^2 \leq \lambda_i(\omh) + \gamma \, \theta^{\frac{1}{2n}},$$
où 
$$B'=B(R')\setminus (B_0 \cap \omh^c), \;\; R'= R +2\theta^{\frac{1}{4n}}.$$ 
\end{lem}   
\begin{rem} Ce lemme est la première étape de la démonstration du fait que le spectre de Dirichlet d'un domaine contenant une boule, dont le volume est presque égal à celui du domaine, est proche du spectre de Dirichlet de cette boule.
\end{rem}
\begin{rem}Un expression explicite des constantes $\gamma$ et $\kappa$ est donnée à la fin de la démonstration.
\end{rem}
\begin{dem}
La construction des fonctions \og tests \fg\, $(F_i)_{1 \leq i \leq k}$ de $H^1_0(B')$ se déroulent en plusieurs étapes. Nous commen\c{c}ons par établir une majoration du terme 
$$\int_{\omh \setminus B(R+\beta)} \phi_i^2$$ 
avec $\phi_i$ une fonction propre de $\omh$ et $\beta$ une constante positive petite (voir (\ref{a4e3})). Nous définissons ensuite des fonctions $(\psi_i)_{1 \leq i \leq k}$. A l'aide de la majoration obtenue dans la première étape, nous estimons l'énergie et la valeur du produit scalaire usuel sur $L^2$ de tout couple de la famille $(\psi_i)_{1 \leq i \leq k}$. Nous terminons la preuve à l'aide d'un lemme technique (lemme \ref{a4l2}).

Dans la démonstration ci-dessous, les constantes $\alpha$ et $\beta$ sont fixées. On les choisira de manière convenable au moment de construire les fonctions tests.\\
Soit $s<t$ deux réels positifs. On note $\chi_{s,t} : \mathbb{R}^n \rightarrow \mathbb{R}^+$ la fonction définie par
$$ \chi_{s,t}(x)= \left\{\begin{array}{rl}
1 & \mbox{ si } |x| \leq s\\
-\frac{|x|}{t-s}+ \frac{t}{t-s} & \mbox{ si } s \leq |x|\leq t \\
0 & \mbox{ si } |x| > t
\end{array}
\right.
$$
Notons $ \chi_1=\chi_{R,R+\beta}$. Soit $\phi_i$ une fonction propre de $\omh$ de valeur propre $\lambda_i(\omh)$ que l'on prolonge par $0$ en dehors de $\omh$. On suppose la fonction $\phi_i$ unitaire pour la norme $L^2$ sur $\omh$. Par construction, $(1-\chi_1)\phi_i$ appartient à $H_0^1(\omh \setminus B_0)$. Pour tout $x$ dans $\omh \setminus B(R+ \beta)$, on a l'égalité
$$ \phi_i^2(x) = ((1-\chi_1)\phi_i)^2(x)$$
d'où 

\begin{equation}\label{a4e4}
\int_{\omh \setminus B(R+ \beta)} \phi_i^2 \leq \int_{\omh} ((1-\chi_1)\phi_i)^2.
\end{equation}

Nous cherchons maintenant à majorer le terme de droite de l'inégalité (\ref{a4e4}). L'inégalité de Faber-Krahn entraine
$$ \lambda_1(\omh \setminus B_0) \geq \frac{c_2(n)}{\theta^{\frac{2}{n}}},$$
avec $c_2(n)= \lambda_1(B(1))w_n^{\frac{2}{n}}$. On en déduit 
$$
\int_{\omh} ((1-\chi_1)\phi_i)^2 \leq \frac{\theta^{\frac{2}{n}}}{c_2(n)}\int_{\omh} |\nabla ((1-\chi_1)\phi_i)|^2.
$$
Or
$$\int_{\omh} |\nabla ((1-\chi_1)\phi_i)|^2 \leq 2\left(\int_{\omh} |\nabla \phi_i|^2 + ||\nabla \chi_1||_{L^{\infty}}^2\int_{\omh \setminus B_0}\phi_i^2\right)$$
et 
$$  ||\nabla \chi_1||_{L^{\infty}}^2= \frac{1}{\beta^2}.$$
On obtient finalement en utilisant (\ref{a4e4})
\begin{equation}\label{a4e3}
\int_{\omh \setminus B(R +\beta)} \phi_i^2 \leq \frac{2\theta^{\frac{2}{n}}}{c_2(n)}\left(\lambda_i(\omh)+ \frac{1}{\beta^2}\right).
\end{equation}

\medskip

Notons $\chi_2= \chi_{R+\beta,R+\beta +\alpha}$ et $\psi_i= \chi_2\phi_i$. Par construction, la fonction $\psi_i$ appartient à $H_0^1(B(R+\alpha+\beta)\setminus (B_0 \cap \omh^c))$. On définit le domaine $B'$ de la manière suivante,  $B'=B(R+\alpha+\beta)\setminus (B_0 \cap \omh^c)$. On a l'égalité
$$ \left|\int_{B'}\psi_i\psi_j-\int_{\omh} \phi_i\phi_j\right|=\left|\int_{\omh} (1-\chi_2^2)\phi_i\phi_j\right|.$$
On déduit de l'inégalité de Cauchy-Schwarz,
$$ \left|\int_{B'}\psi_i\psi_j-\int_{\omh} \phi_i\phi_j\right| \leq \left(\int_{\omh \setminus B(R+\beta)}\phi_i^2\right)^{\moi}\left(\int_{\omh \setminus B(R+\beta)}\phi_j^2\right)^{\moi}.$$

Estimons maintenant l'énergie des fonctions $\psi_i$.
\begin{multline*}
 \int_{B'}|\nabla \psi_i|^2 \leq \lambda_i(\omh) + 
 ||\nabla \chi_2||^2_{\infty}\int_{\omh \setminus B(R+\beta)}\phi_i^2 + \\
 2 ||\nabla \chi_2||_{\infty}\sqrt{\lambda_i(\omh)} \left(\int_{\omh \setminus B(R+\beta)}\phi_i^2\right)^{\moi}.
 \end{multline*}
On déduit de (\ref{a4e3}) 
\begin{multline*}
 \int_{B'}|\nabla \psi_i|^2 \leq \lambda_i(\omh) + \frac{2\theta^{\frac{2}{n}}}{\alpha^2 c_2(n)}\left(\lambda_k(\omh)+\frac{1}{\beta^2}\right)+ \\
 \frac{2}{\alpha}\left(\lambda_k(\omh)\frac{2\theta^{\frac{2}{n}}}{c_2(n)}\left(\lambda_k(\omh)+ \frac{1}{\beta^2}\right)
\right)^{\moi}
\end{multline*} 
et
$$ \left|\int_{B'}\psi_i\psi_j-\int_{\omh} \phi_i\phi_j\right| \leq \left(\lambda_k(\omh)+\frac{1}{\beta^2}\right)\frac{2\theta^{\frac{2}{n}}}{c_2(n)}.$$
Fixons $\alpha^2=\beta^2=\theta^{\frac{1}{2n}}$. On en déduit (on suppose $\theta <1$)
$$ \int_{B'}|\nabla \psi_i|^2 \leq \lambda_i(\omh) + \frac{2}{c_2(n)}\big(\lambda_k(\omh)+1\big)\theta^{\frac{1}{n}} + 2\left(\lambda_k(\omh)\frac{2}{c_2(n)}\left(\lambda_k(\omh)+ 1\right)\right)^{\moi}\, \theta^{\frac{1}{2n}}.$$
On obtient
$$ \int_{B'}|\nabla \psi_i|^2 \leq \lambda_i(\omh) +c(k,n)\,\theta^{\frac{1}{2n}},$$
avec 
\begin{equation}\label{a10e1}
c(k,n)=\left(\frac{8 }{c_2(n)}+1\right)\max\{\lambda_k(\omh),1\}.
\end{equation}
On en déduit également
$$ \left|\int_{B'}\psi_i\psi_j-\int_{\omh} \phi_i\phi_j\right| \leq c(k,n)\,\theta^{\frac{1}{2n}}.$$

\smallskip

On termine la preuve du lemme \ref{a4l3} à l'aide du lemme suivant, dont une démonstration se trouve en annexe. 
\begin{lem}\label{a4l2}Soit $(E,\langle\cdot,\cdot\rangle )$ un espace euclidien de dimension $k$ et $q$ une forme quadratique sur $E$. Soit $(f_i)_{1\leq i \leq k}$ une famille d'éléments de $E$ vérifiant pour tout $i,j$ dans $\{1,\cdots,k\}$ 
$$ |\langle f_i,f_j\rangle -\delta_{i,j}| \leq c \mbox{ et } q(f_i) \leq \lambda_i + c$$
où $\delta_{i,j}$ désigne le symbole de Kronecker, $c$ un réel positif et $\lambda_i$ une famille croissante de nombres positifs. Supposons de plus que le réel $c$ vérifie $ca_k\leq \frac{1}{4}$ où $a_k$ est défini par la relation de récurrence $a_s=1+\sum_{i=1}^{s-1}a_i^2$ et $a_1=1$. Sous ces hypothèses, il existe une base orthonormée $(F_i)_{1\leq i \leq k}$ de $(E, \langle\cdot,\cdot\rangle )$ telle que pour tout $i$ dans $\{1,\cdots,k\}$,
$$ q(F_i) \leq \lambda_i + 14ka_k\max\{\lambda_k,1\}c.$$
\end{lem}
En appliquant le lemme \ref{a4l2} aux fonctions $(\psi_i)$, on obtient une famille orhonormée de fonctions $(F_i)_{1 \leq i \leq k}$ qui vérifient pour $\theta < \kappa$ 
$$
q(F_i) \leq \lambda_i(\om)+ \gamma \, \theta^{\frac{1}{2n}},
$$
avec 
$$\kappa=\kappa(k,n,\lambda_k(\omh))= \left(\frac{1}{4a_k c(k,n) }\right)^{2n}$$
où $c(k,n)$ est la constante définie par (\ref{a10e1}) et 
$$\gamma=\gamma (k,n,\lambda_k(\omh))= \left(\frac{8}{c_2(n)}+1\right)14ka_k\max\{\lambda_k(\omh)^2,1\},$$
où $c_2(n)= \lambda_1(B(1))w_n^{\frac{2}{n}}$.

\end{dem}

En appliquant le lemme \ref{a4l3} dans le cas où $\omh=\om$ et en utilisant l'estimation de volume (\ref{a4e10}), on en déduit le
\begin{cor}\label{a4c1} Soit $\om$ un domaine vérifiant les hypothèses du théorème \ref{convergence}. En conservant les notations du lemme \ref{a4l3}, on a sous l'hypothèse \\$\ep <\min\{1,\rho(vol (\om),n,\lambda_1(B_0))\}$, l'inégalité
\begin{equation}\label{a4e16}
\lambda_1(B(R')\setminus (\om^c \cap B_0)) \leq \lambda_1(\om) + \gamma(1,n,2\lambda_1(B_0))\big(c_1(n)vol(\om)\big)^{\frac{1}{2n}}\ep^{\frac{1}{20n}},
\end{equation}
avec
$$ R'=R + 2 \big(c_1(n)vol (\om)\big)^{\frac{1}{4n}}\ep^{\frac{1}{40n}} \mbox{ et } \rho(vol (\om),n,\lambda_1(B_0))=\left(\frac{\kappa(1,n,2\lambda_1(B_0))}{c_1(n)vol (\om)}\right)^{10}.$$
Pour alléger les notations, nous noterons
$$\tau(\ep) =\gamma(1,n,2\lambda_1(B_0))\big(c_1(n)vol(\om)\big)^{\frac{1}{2n}}\ep^{\frac{1}{20n}}.$$

\end{cor}

\subsection{Majoration de l'expression $\lambda_k(B_0)- \lambda_k(\om \cup B_0)$}
Le lemme \ref{a4l3} permet également de démontrer le résultat suivant.
\begin{lem}\label{lemme3.13}Soit $\om$ un domaine vérifiant les hypothèses du théorème \ref{convergence}. Nous conservons les notations du lemme \ref{a4l3}. Pour tout entier $k \geq 2$ et pour $\ep < \mu (\lambda_k(B_0),k,n,vol (\om))$, il existe une famille de fonctions $(F_i)_{1 \leq i \leq k}$ de $H^1_0(B_0)$, orthonormée pour le produit scalaire $L^2$ sur $B_0$, qui vérifie pour tout $i$ dans $\{1,\cdots,k\}$, 
$$ \int_{B_0} |\nabla F_i|^2 \leq \lambda_i(\om \cup B_0) + \beta \big(\lambda_k(B_0),k,n,vol (\om)\big) \ep^{\frac{1}{40n}},$$
avec 
$$\beta\big(\lambda_k(B_0),k,n, vol (\om)\big)= c_4(k,n)\max\{1,\lambda_k^2(B_0)\}\max \{vol(\om)^{-\frac{3}{2n}},vol(\om)^{\frac{1}{2n}}\}$$
 et 
$$\mu(\lambda_k(B_0),k,n,vol (\om))=\left(\frac{\kappa(k,n,\lambda_k(B_0))}{c_1(n)vol (\om)}\right)^{10}.$$
\end{lem}
\begin{dem}
Il suffit d'appliquer une homothétie aux fonctions cons\-truites dans le lemme \ref{a4l3} dans le cas où $\omh=\om \cup B_0$ et de remarquer que $\lambda_k(\om \cup B_0) \leq \lambda_k(B_0)$.
\end{dem}

Nous déduisons du lemme \ref{lemme3.13}, une majoration de $\lambda_k(B_0)- \lambda_k(\om \cup B_0)$ à l'aide du lemme \ref{a4l1} dont une preuve se trouve en annexe. 
\begin{defi}Soit $(\mu_i)_{i \geq 1}$ une suite croissante de nombres réels positifs. Pour tout entier $i$ positif, on note
$$\mu_i^+=\min \{\mu_j;\;\;\mu_j>\mu_i\}.$$
\end{defi}

\begin{lem}\label{a4l1}
Soit $(H,\langle\cdot,\cdot\rangle )$ un espace de Hilbert et $q$ une forme quadratique définie sur un domaine $D(q)$ dense dans $H$. On suppose que le spectre de $q$ est discret, on note $(\mu_i)_{i \geq 1}$ les valeurs propres de $q$ et $(h_i)_{i \geq 1}$ une base orthonormée de vecteurs propres de $q$. Soit $k$ un entier positif fixé. Supposons qu'il existe une famille $(f_i)_{1 \leq i \leq k}$ orthonormée qui vérifie pour tout $i$ dans $\{1,\cdots,k\}$, 
 $$q(f_i) \leq \lambda_i + \eta$$
avec $0<\lambda_i \leq \mu_i$. Il existe une constante $c_k$ telle que si $\eta$ vérifie $\eta \leq \frac{1}{2c_k}$ alors on a  pour tout $i$ dans $\{1,\cdots,k\}$,
$$ \mu_i \leq \lambda_i + c_{k+1}\eta.$$
La constante $c_k$ est définie par $c_k=(8t_k\mu_k)^k\frac{1}{\mu_k}$ où \\$t_k= \max (\max_{1\leq i \leq k}(\frac{1}{\mu_i^+ - \mu_i}),1)$.
\end{lem}

En appliquant ce lemme avec $\lambda_i=\lambda_i(\om \cup B_0)$ et $\mu_i=\lambda_i(B_0)$, on en déduit la proposition suivante
\begin{prop}\label{a4p1}
Soit $\om$ un domaine vérifiant les hypothèses du théorème \ref{convergence}. Supposons que $\epsilon <\min\left\{\mu,\left(\frac{1}{2c_k\beta}\right)^{40n}\right\}$  alors pour tout $i$ dans $\{1,\cdots,k\}$, on a 
$$ \lambda_i(B_0) \leq \lambda_i(B_0 \cup \om) + c_{k+1}\beta\,\ep^{\frac{1}{40n}}$$
où $c_k$ est la constante définie dans le lemme \ref{a4l1}, $\mu=\mu(\lambda_k(B_0),k,n,vol (\om))$, $\beta= \beta\big(\lambda_k(B_0),k,n, vol (\om)\big)$ sont les constantes définies dans le lemme \ref{lemme3.13}.
\end{prop}

\subsection{Majoration de l'expression $\lambda_k( D_1)-\lambda_k(B_0)$}

Commen\c{c}ons par faire la remarque suivante (en utilisant le théorème \ref{capacite} dans le cas où $\om=B(R')$ et les notations introduites dans son énoncé).
$$
\begin{array}{rcl}
\lambda_k(D_1)-\lambda_k(B_0) & \leq & \lambda_k(D_1)-\lambda_k(B(R')) \\
                              & \leq & C_k \big(\ca (A_{R,R'} \cup A_0)\big)^{\moi}
\end{array}
$$                              
si $\ca (A_{R,R'} \cup A_0)$ est assez petit (le domaine $A_{R,R'}$ désigne $B(R')\setminus B(R)$ et $R'$ est défini dans le corollaire \ref{a4c1}). La capacité de Dirichlet est sous-additive (lemme \ref{a4l6}), il suffit donc de majorer séparément les deux termes. Nous utilisons de nouveau le théorème \ref{capacite}.
$$\ca (A_{R,R'}) \leq B_1\lambda_1(B_0)\left(1-\left(\frac{R}{R'}\right)^2\right)$$
et
$$\ca(A_0) \leq B_1\big(\lambda_1(B(R')\setminus A_0)- \lambda_1(B(R'))\big).$$
D'où en utilisant le corollaire \ref{a4c1},
$$\ca(A_0) \leq B_1\left(\ep\lambda_1(B_0) + \tau(\ep) + \lambda_1(B_0)\left(1-\left(\frac{R}{R'}\right)^2\right)\right),$$
par hypothèse sur la première valeur propre $\lambda_1(\om)$.\\
Notons $c_5(n,vol (\om))=2\big(c_1(n)vol (\om)\big)^{\frac{1}{4n}}$, nous avons l'égalité
$$ \frac{R'}{R}= 1 + \frac{1}{R}c_5(n, vol (\om))\ep^{\frac{1}{40n}},$$
d'où
$$ \left(\frac{R}{R'}\right)^2 \geq \left(1 -\frac{1}{R}c_5(n, vol (\om))\ep^{\frac{1}{40n}}\right)^2,$$   
on en déduit
$$ 1-\left(\frac{R}{R'}\right)^2 \leq 4\left(\frac{c_1(n)}{w_n^4}\right)^{\frac{1}{4n}}vol (\om)^{-\frac{3}{4n}}\ep^{\frac{1}{40n}}.$$
Finalement, il existe une constante $\delta (n,vol (\om), \lambda_1(B_0))$ telle que
$$ \ca (A_0) \leq \delta (n,vol (\om), \lambda_1(B_0))\ep^{\frac{1}{40n}}.$$
Par conséquent si $\delta (n,vol (\om), \lambda_1(B_0))\ep^{\frac{1}{40n}}\leq \ep_k$, on a 
$$\lambda_k(D_1)-\lambda_k(B_0) \leq C_k\delta^{\moi} (n,vol (\om), \lambda_1(B_0))\ep^{\frac{1}{80n}}.$$
 
\subsection*{Annexe : Démonstration des lemmes \ref{a4l2} et \ref{a4l1}}
\subsubsection*{Démonstration du lemme \ref{a4l2}}
Rappelons l'énoncé du lemme \ref{a4l2}.
\begin{lembis}Soit $(E,\langle\cdot,\cdot\rangle )$ un espace euclidien de dimension $k$ et $q$ une forme quadratique sur $E$. Soit $(f_i)_{1\leq i \leq k}$ une famille d'éléments de $E$ vérifiant pour tout $i,j$ dans $\{1,\cdots,k\}$ 
$$ |\langle f_i,f_j\rangle -\delta_{i,j}| \leq c \mbox{ et } q(f_i) \leq \lambda_i + c$$
où $\delta_{i,j}$ désigne le symbole de Kronecker, $c$ un réel positif et $\lambda_i$ une famille croissante de nombres positifs. Supposons de plus que le réel $c$ vérifie $ca_k\leq \frac{1}{4}$ où $a_k$ est défini par la relation de récurrence $a_s=1+\sum_{i=1}^{s-1}a_i^2$ et $a_1=1$. Sous ces hypothèses, il existe une base orthonormée $(F_i)_{1\leq i \leq k}$ de $(E, \langle\cdot,\cdot\rangle )$ telle que pour tout $i$ dans $\{1,\cdots,k\}$,
$$ q(F_i) \leq \lambda_i + 14ka_k\max\{\lambda_k,1\}c.$$
\end{lembis}
\begin{dem}
Dans la démonstration, nous utiliserons les majorations  ci-dessous, valables sous l'hypothèse $0 \leq c \leq \moi$.
$$\frac{1}{1-c}\leq 1+4c, \; \sqrt{\frac{1}{1-c}}\leq 1+2c.$$
Notons
$$ F_1=\frac{f_1}{||f_1||},$$
où $|| \cdot ||$ est la norme associée au produit scalaire $\langle\cdot,\cdot\rangle $.\\
Pour tout $i$ dans $\{2,\cdots,k\}$, on définit par récurrence
$$h_i = f_i-\sum_{j=1}^{i-1}\langle F_j,f_i\rangle F_j
\mbox{  et  }  F_i=\frac{h_i}{||h_i||}.$$
Pour tout $s$ dans $\{2,\cdots,k\}$ et sous l'hypothèse $c \leq \moi$, on a
$$ |\langle f_s,F_1\rangle | \leq \sqrt{2}c.$$
Nous allons montrer par récurrence que pour tout $i$ dans $\{1,\cdots,k-1\}$ et pour tout $s>i$, on a
$$ |\langle f_s,F_i\rangle | \leq \sqrt{2}a_ic.$$
Fixons $i$ et $s$ tels que $i<s$, alors les inégalités ci-dessous sont vérifiées
$$ |\langle f_s,F_i\rangle | \leq \frac{1}{||h_i||}\left(|\langle f_s,f_i\rangle | + \sum_{j=1}^{i-1}|\langle f_s,F_j\rangle \langle f_i,F_j\rangle |\right)$$
et 
$$ ||h_i|| \geq \left(1-c-\sum_{j=1}^{i-1}|\langle f_s,F_j\rangle \langle f_i,F_j\rangle |\right)^{\moi}.$$
Par conséquent, par hypothèse de récurrence et sous l'hypothèse $a_{i-1}c\leq \moi$, 
$$|\langle f_s,F_i\rangle |\leq \frac{1}{||h_i||}\Big(c + 2c^2 \sum_{j=1}^{i-1} a_j^2\Big).$$
Notons $p(f_i)= \sum_{j=1}^{i-1}\langle F_j,f_i\rangle F_j$. On déduit de l'égalité $\langle h_i,p(f_i)\rangle =0$
$$||h_i||^2= ||f_i||^2 - ||p(f_i)||^2,$$
d'où
$$ ||h_i|| \geq \Big(1-c -2 c^2 \sum_{j=1}^{i-1} a_j^2\Big)^{\moi}.$$
Par conséquent, si $a_ic \leq \moi$ on en déduit (en utilisant $c \leq \moi$)
$$ |\langle f_s,F_i\rangle |\leq \sqrt{2}c\Big(1+ \sum_{j=1}^{i-1}a_j^2\Big),$$
ce qui démontre la propriété.

Estimons maintenant $q(F_i)$ pour $i$ appartenant à $\{1,\cdots,k\}$. Par définition de $F_1$, on a
$$q(F_1)\leq \frac{\lambda_1}{1-c} +\frac{c}{1-c},$$
d'où
$$ q(F_1) \leq \lambda_1 + 4\lambda_kc +2c.$$
Nous allons montrer par récurrence que l'hypothèse $ca_k \leq \frac{1}{4}$ entraîne 
$$ q(F_k)-\lambda_k \leq 14cka_k\max\{\lambda_k,1\}.$$
L'inégalité est vérifiée pour $k=1$. Pour alléger les notations, notons
$$A_k= \lambda_k + 14cka_k\max\{\lambda_k,1\}.$$
 Fixons $s$ tel que $2\leq s\leq k$. Par définition de $F_s$, on a 
$$q(F_s) \leq \frac{1}{||h_s||^2}\Big(q(f_s)+q(p(f_s))+2\sqrt{q(p(f_s))q(f_s)}\Big).$$
En utilisant la forme bilinéaire associée à $q$, on obtient
$$ q(p(f_s)) \leq A_{s-1}\sum_{1\leq i,j \leq s-1}|\langle f_s,F_i\rangle \langle f_s,F_j\rangle |.$$
En utilisant l'inégalité de Cauchy-Schwarz, il vient
$$ q(p(f_s)) \leq 2A_{s-1}c^2(a_s-1).$$
On déduit de l'hypothèse $q(f_s) \leq \lambda_s +c$, de $||h_s||^2 \geq 1-a_sc$ et de $a_sc \leq \moi$,
$$ q(F_s) \leq \lambda_s(1+4a_sc) +2\big(c + 2c^2A_{s-1}(a_s-1)+2c\sqrt{2A_{s-1}(a_s-1)}\sqrt{q(f_s)}\big),$$
que l'on majore par
$$ q(F_s)-\lambda_s \leq 4a_sc\lambda_k +2\Big(c +cA_{s-1} +c\big(2A_{s-1}(a_s-1) + \lambda_s +c\big)\Big).$$
On en déduit 
$$ q(F_s) -\lambda_s \leq 10ca_s\max\{\lambda_k,1\} +4ca_sA_{s-1},$$
on conclut en utilisant $4ca_s \leq 1$.
\end{dem}

\subsubsection*{Démonstration du lemme \ref{a4l1}} 
Démontrons le
\begin{lembis}Soit $(H,\langle\cdot,\cdot\rangle )$ un espace de Hilbert et $q$ une forme quadratique définie sur un domaine $D(q)$ dense dans $H$. On suppose que le spectre de $q$ est discret, on note $(\mu_i)_{i \geq 1}$ les valeurs propres de $q$ et $(h_i)_{i \geq 1}$ une base orthonormée de vecteurs propres de $q$. Soit $k$ un entier positif fixé. Supposons qu'il existe une famille $(f_i)_{1 \leq i \leq k}$ orthonormée qui vérifie pour tout $i$ dans $\{1,\cdots,k\}$, 
 $$q(f_i) \leq \lambda_i + \eta$$
avec $0<\lambda_i \leq \mu_i$. Il existe une constante $c_k$ telle que si $\eta$ vérifie $\eta \leq \frac{1}{2c_k}$ alors on a  pour tout $i$ dans $\{1,\cdots,k\}$,
$$ \mu_i \leq \lambda_i + c_{k+1}\eta.$$
La constante $c_k$ est définie par $c_k=(8t_k\mu_k)^k\frac{1}{\mu_k}$\\ où $t_k= \max \big(\max_{1\leq i \leq k}\big(\frac{1}{\mu_i^+ - \mu_i}\big),1\big)$.
\end{lembis}
\begin{dem}

Pour tout $i$ dans $\{1,\cdots,k\}$, on décompose la fonction $f_i$ comme suit
$$ f_i= \sum_{s=1}^{i}a_{is}h_s + b_ig_i,$$
où l'on peut supposer que la fonction $g_i$ est de norme $1$ et $q-$orthogonal aux $i$ premiers espaces propres de $q$. Sous ces hypothèses, on a donc en particulier 
\begin{equation}\label{a4e14}
\sum_{s=1}^{i}a_{is}^2 + b_i^2=1.
\end{equation}
Pour démontrer le lemme, il suffit de montrer que les coefficients $a_{ii}$ sont proches de $1$.

L'orthonormalité de la famille $(f_i)_{1 \leq i \leq k}$ implique (pour $j<i$) l'estimation
\begin{equation}\label{a4e12}
 |a_{ij}a_{jj}| \leq \left(\sum_{s=1}^{j-1}a_{is}^2\right)^{\moi}\left(\sum_{s=1}^{j-1}a_{js}^2\right)^{\moi}+|b_j|.
\end{equation} 
En particulier on a pour tout $i$ dans  $\{2,\cdots,k\}$,
\begin{equation}\label{a4e5}
a_{i1}^2\leq \frac{b_1^2}{1-b_1^2}.
\end{equation}
Notons pour $i>j$,
$$C_{i,j}=\sum_{s=1}^{j}a_{is}^2.$$
Montrons par récurrence que pour tout $i>j$,
\begin{equation}\label{a4e6}
 C_{i,j} \leq D_j,
\end{equation} 
avec $D_j=  (8t_k\mu_k)^j\frac{1}{\mu_k}\eta$.

Estimons maintenant, les coefficients $b_j$.
$$ q(f_i)\geq \sum_{s=1}^i a_{is}^2\mu_s +b_i^2\mu_i^+.$$
on réécrit cette inégalité sous la forme 
$$ q(f_i)-\mu_i \geq \sum_{s=1}^{i-1} a_{is}^2(\mu_s-\mu_i) +b_i^2(\mu_i^+ -\mu_i).$$
On en déduit
\begin{equation}\label{a4e11}
 b_1^2 \leq t_k\eta
\end{equation}
et pour $i>1$,
$$
 b_i^2 \leq \frac{1}{\mu_i^+ -\mu_i}(C_{i,i-1}\mu_i+ \eta),
$$
ce qui implique
\begin{equation}\label{a4e7}
 b_i^2 \leq t_k\mu_k C_{i,i-1}+ t_k\eta.
\end{equation}
On déduit de (\ref{a4e11}) et (\ref{a4e5}) que pour tout $i$ dans $\{2,\cdots,k\}$, $C_{i,1} \leq D_1$. En effet par définition de $t_k$, $\mu_kt_k \geq 1$ et donc $t_k\eta \leq c_1\eta \leq c_k\eta \leq \moi$. 
Fixons $i$ et $j$ tels que $i>j$.
En utilisant l'hypothèse de récurrence et (\ref{a4e12}), on obtient
$$
  C_{i,j} \leq D_{j-1}+ \frac{2}{a_{jj}^2}(D_{j-1}^2+ b_j^2).
$$
En utilisant (\ref{a4e7}), il vient
$$ C_{i,j} \leq D_{j-1} + \frac{2}{a_{jj}^2}D_{j-1}^2+ \frac{2}{a_{jj}^2}(t_k\mu_kD_{j-1}+ t_k\eta).$$
On en déduit en utilisant $D_j \leq c_k\eta \leq \moi$,
$$ C_{i,j} \leq  D_{j-1}\left(1 + \frac{1}{a_{jj}^2}+ \frac{2}{a_{jj}^2}t_k\mu_k\right) + \frac{2}{a_{jj}^2}t_k\eta.$$
On obtient également en utilisant (\ref{a4e14}) et (\ref{a4e7}),
$$ a_{jj}^2 \geq 1 - \big(D_{j-1}(1+t_k\mu_k)+t_k\eta\big).$$
En remarquant que $D_{j-1}(1+t_k\mu_k)+t_k\eta \leq D_j\leq \moi$, on en déduit
\begin{equation}\label{a4e13}
a_{jj}^2 \geq 1 - D_j
\end{equation}
et
$$C_{i,j} \leq D_{j-1}(3 + 4t_k\mu_k) +4t_k\eta,$$
d'où en majorant grossièrement,
$$C_{i,j} \leq 8t_k\mu_kD_{j-1},$$
ce qui démontre (\ref{a4e6}). Pour terminer la preuve du lemme, on remarque que
$$ a_{jj}^2\mu_j \leq q(f_j),$$
d'où l'on tire par (\ref{a4e13})
$$ \mu_j \leq \frac{\lambda_j}{1-D_j} +2\eta,$$
en utilisant de nouveau $D_j \leq \moi$ et $\lambda_j \leq \mu_j$, on obtient
$$ \mu_j \leq \lambda_j + 4D_j\mu_k +2 \eta,$$
on conclut en majorant convenablement.
\end{dem}

\bibliographystyle{plain}

\medskip

\noindent
J. Bertrand\\
{\it Institut Fourier, Université Joseph Fourier, Grenoble, France.}\\
Adresse électronique : jbertran@fourier.ujf-grenoble.fr

\medskip

\noindent
B. Colbois\\
{\it Institut de Math\'ematiques, Universit\'e de Neuch\^atel, Neuchâtel, Suisse.}\\
Adresse électronique : bruno.colbois@unine.ch

\end{document}